\documentclass[manuscript,screen,nonacm]{acmart}

\AtBeginDocument{%
  \providecommand\BibTeX{{%
    \normalfont B\kern-0.5em{\scshape i\kern-0.25em b}\kern-0.8em\TeX}}}

\setcopyright{none}
\copyrightyear{2021}
\acmYear{2021}
\acmDOI{}

\acmConference[Submitted]{}{}{}
\acmBooktitle{}
\acmPrice{}
\acmISBN{}

\begin{document}

\title{The Old and the New: Can Physics-Informed Deep-Learning Replace Traditional Linear Solvers?}

\author{Stefano Markidis}
\email{markidis@kth.se}
\affiliation{%
  \institution{KTH Royal Institute of Technology}
  \city{Stockholm}
  \country{Sweden}
}

\renewcommand{\shortauthors}{Markidis}

\begin{abstract}
Physics-Informed Neural Networks (PINN) are neural networks encoding the problem governing equations, such as Partial Differential Equations (PDE), as a part of the neural network. PINNs have emerged as a new essential tool to solve various challenging problems, including computing linear systems arising from PDEs, a task for which several traditional methods exist. In this work, we focus first on evaluating the potential of PINNs as linear solvers in the case of the Poisson equation, an omnipresent equation in scientific computing. We characterize PINN linear solvers in terms of accuracy and performance under different network configurations (depth, activation functions, input data set distribution). We highlight the critical role of transfer learning. Our results show that low-frequency components of the solution converge quickly as an effect of the F-principle. In contrast, an accurate solution of the high frequencies requires an exceedingly long time. To address this limitation, we propose integrating PINNs into traditional linear solvers. We show that this integration leads to the development of new solvers whose performance is on par with other high-performance solvers, such as PETSc conjugate gradient linear solvers, in terms of performance and accuracy. Overall, while the accuracy and computational performance are still a limiting factor for the direct use of PINN linear solvers, hybrid strategies combining old traditional linear solver approaches with new emerging deep-learning techniques are among the most promising methods for developing a new class of linear solvers. 
\end{abstract}

\maketitle

\section{Introduction}
Deep Learning (DL) has revolutionized the way of performing classification, pattern recognition, and regression tasks in various application areas, such as image and speech recognition, recommendation systems, natural language processing, drug discovery, medical imaging,  bioinformatics, and fraud detection, among few examples~\cite{goodfellow2016deep}. However, scientific applications solving linear and non-linear equations with demanding accuracy and computational performance requirements have not been the DL focus. Only until recently, a new class of DL networks, called \emph{Physics-Informed Neural Networks} (PINN), emerged as a very promising DL method to solve scientific computing problems~\cite{raissi2019physics, raissi2017physicsI, raissi2017physicsII}. In fact, PINNs are specifically designed to integrate scientific computing equations, such as Ordinary Differential Equations (ODE), Partial Differential Equations (PDE), non-linear and integral-differential equations~\cite{pang2019fpinns}, into the DL network training. In this work, we focus on PINN application to solve a traditional scientific computing problem: the solution of a linear system arising from the discretization of a PDE. We solve the linear system arising from the Poisson equation, one of the most common PDEs whose solution still requires a non-negligible time with traditional approaches. We evaluate the level of maturity in terms of accuracy and performance of PINN linear solver, either as a replacement of other traditional scientific approaches or to be deployed in combination with conventional scientific methods, such as the multigrid and Gauss-Seidel methods~\cite{quarteroni2010numerical}.

PINNs are deep-learning networks that, after training (solving an optimization problem to minimize a residual function), output an approximated solution of differential equation/equations, given an input point in the integration domain (called collocation point). Before PINNs, previous efforts, have explored solving PDEs with constrained neural networks~\cite{lagaris1998artificial,psichogios1992hybrid}. The major innovation with PINN is the introduction of a \emph{residual} network that encodes the governing physics equations, takes the output of a deep-learning network (called \emph{surrogate}), and calculates a residual value (a loss function in DL terminology). The inclusion of a \emph{residual} network, somehow, bears a resemblance of those iterative Krylov linear solvers in scientific applications. The fundamental difference is that PINNs calculate differential operators on graphs using automatic differentiation~\cite{baydin2018automatic} while traditional scientific approaches are based on numerical schemes for differentiation. As noted in previous works~\cite{raissi2019physics,mishra1}, automatic differentiation is the main strength of PINNs because operators on the residual network can be elegantly and efficiently formulated with automatic differentiation. An important point is that the PINN's \emph{residual} network should not be confused with the popular network architectures, called also \emph{Residual} networks, or \emph{ResNet} in short, where the name derives from using skip-connection or residual connections~\cite{goodfellow2016deep} instead of calculating a residual like in PINNs.

\textbf{The basic formulation of the PINN training does not require labeled data, e.g., results from other simulations or experimental data, and is unsupervised}: PINNs only require the evaluation of the residual function~\cite{mishra1}. Providing simulation or experimental data for training the network in a supervised manner is also possible and necessary for so data-assimilation~\cite{raissi2020hidden}, inverse problems~\cite{mishra2}, super resolution~\cite{esmaeilzadeh2020meshfreeflownet,wang2020physics}, and discrete PINNs~\cite{raissi2019physics}. The supervised approach is often used for solving ill-defined problems when for instance we lack boundary conditions or an Equation of State (EoS) to close a system of equations (for instance, EoS for the fluid equations~\cite{zhu2020generating}). In this study, we only focus on the basic PINNs as we are interested in solving PDEs without relying on other simulations to assist the DL network training. A common case in scientific applications is that we solve the same PDE with different source terms at each time step. For instance, in addition to other computational kernels, Molecular Dynamics (MD) code and semi-implicit fluid and plasma codes, such as GROMACS~\cite{van2005gromacs}, Nek5000~\cite{nek5000-web-page}, and iPIC3D~\cite{markidis2010multi}, calculate the Poisson equation for the electrostatic and pressure solver~\cite{offermans2016strong} and divergence cleaning operations at each cycle.

Once a PINN is trained, the inference from the trained PINN can be used to replace traditional numerical solvers in scientific computing. In this so-called \emph{inference} or \emph{prediction} step, the input includes independent variables like simulation time step and simulation domain positions. The output is the solution of the governing equations at the time and position specified by the input. Therefore, PINNs are a \emph{gridless} method because any point in in the domain can be taken as input without requiring the definition of a mesh. Moreover, the trained PINN network can be used for predicting the values on simulation grids of different resolutions without the need of being retrained. For this reason, the computational cost does not scale with the number of grid points like many traditional computational methods. PINNs  borrow concepts from popular methods in traditional scientific computing, including Newton-Krylov solvers~\cite{kelley1995iterative}, finite element methods (FEM)~\cite{rao2017finite}, and Monte Carlo techniques~\cite{rubinstein2016simulation}. Like the Newton-Krylov solvers, PINNs training is driven by the objective of minimizing the residual function and employs Newton methods during the optimization process. Similarly to the FEM, PINN uses interpolation basis (non-linear) functions, called \emph{activation functions}~\cite{ramachandran2017searching} in the neural network fields. Like Monte Carlo and quasi-Monte Carlo methods, PINNs integrate the governing equations using a random or a low-discrepancy sequence, such as the Sobol sequence~\cite{sobol1990quasi}, for the collocation points used during the evaluation the residual function. 

The motivation of this work is twofold. First, we evaluate the potential of deploying PINNs for solving linear systems, such as the one arising from the Poisson equation. We focus on solving the Poisson equation, a generalization of the Laplace equation, and an omnipresent equation in scientific computing. Traditionally, Poisson solvers are based on linear solvers, such as the Conjugate Gradient (CG) or Fast Fourier Transform (FFT). These approaches may require a large number of iterations before convergence and are computationally expensive as the fastest methods scale as $\mathcal{O}(N_g \log N_g)$, where $N_g$ is the number of grid points in the simulation domain. The second goal of this work is to propose a new class of linear solvers combining new emerging DL approaches with old traditional linear solvers, such as multigrid and iterative solvers.

In this work, we show that the accuracy and the convergence of PINN solvers can be tuned by setting up an appropriate configuration of depth, layer size, activation functions and by leveraging transfer learning. We find that fully-connected surrogate/approximator networks with more than three layers produce similar performance results in the first thousand training epochs. The choice of activation function is critical for PINN performance: depending on the \emph{smoothness} of the source term, different activation functions provide considerably different accuracy and convergence. Transfer learning in PINNs allow us to initialize the network with the results of another training solving the same PDE with a different source term~\cite{weiss2016survey}. The usage of transfer learning considerably speed-up the training of the network. In terms of accuracy and computational performance, a naive replacement of traditional numerical approaches with the direct usage of PINNs is still not competitive with traditional solvers and codes, such as CG implementations in HPC packages~\cite{balay2019petsc}. 

To address the limitations of the direct usage of PINN, we combine PINN linear solvers with traditional approaches such as the multigrid and Gauss-Seidel methods~\cite{trottenberg2000multigrid,quarteroni2010numerical}. The DL linear solver is used to solve the linear system on a coarse grid and the solution refined on finer grids using the multigrid V-cycle and Gauss-Seidel solver iterations. This approach allows us to use the DL networking of converging quickly on low-frequency components of the problem solution and rely on Gauss-Seidel to solve accurately high-frequency components of the solution. We show that the integration of DL techniques in traditional linear solvers leads to solvers that are on-par of high-performance solvers, such as PETSc conjugate gradient linear solvers, both in terms of performance and accuracy.

The paper is organized as follows. We first introduce the governing equations, the background information about PINN architecture and showcase the usage of PINN to solve the 2D Poisson equation. Section~\ref{sec:tune} presents a characterization of PINN linear solver performance when varying the network size, activation functions, and data set distribution and we highlight the critical importance of leveraging transfer learning. We present the design of a Poisson solver combining new emerging DL techniques into the V-cycle of the multigrid method and analyze its error and computational performance in Section~\ref{sec:integrate}. Finally, we summarize this study and outline challenges and next step for the future work in Section~\ref{sec:conclusion}.

\section{The New: Physics-Informed Linear Solvers}\label{sec:bg}
The PINNs goal is to approximate the solution of a system of one or more differential, possibly non-linear equations, by encoding explicitly the differential equation formulation in the neural network. Without loss of generality, PINN solves the non-linear equation:
\begin{equation}
u(x)_t = \mathcal{N}u(x)  = 0, x \in \Omega, t \in [0, T],
\end{equation}
where $u$ is the solution of the system, $u_t$ is its derivative with respect to time $t$ in the period [0, T], $\mathcal{N}$ is a non-linear differential operator, $x$ is an independent, possibly multi-dimensional variable, defined over the domain $\Omega$. As a main reference equation to solve, we consider the Poisson equation in a unit square domain and Dirichlet boundary conditions throughout this paper:
\begin{equation}
\nabla^2 u(x,y) = f(x,y), (x,y) \in [0, 1] \times [0, 1] .
\label{poisson}
\end{equation}
While this problem is linear in nature and PINNs can handle non-linear problems, we focus on the Poisson equation because it is one of the most solved PDEs in scientific applications. The Poisson equation, an example of elliptic PDE, arises in several different fields from electrostatic problems in plasma and MD codes, to potential flow and pressure solvers in Computational Fluid Dynamics (CFD), to structural mechanics problems. Elliptic problems are one of the Achilles' heels for scientific applications~\cite{morton2005numerical}. While relatively fast and straightforward - albeit subject to numerical constraints - computational methods exist for solving hyperbolic and parabolic problems, e.g. explicit differentiation, traditionally the solution of elliptic problems requires linear solvers, such as Krylov (CG or GMREs) solvers or FFT. Typically, in scientific applications, the simulation progresses through several time steps, where a Poisson equation with same boundary conditions and different source term $f(x,y)$ (typically not considerably different from the source term of the previous time step) is solved.

\begin{figure}[h!]
\begin{center}
\includegraphics[width=0.8\textwidth]{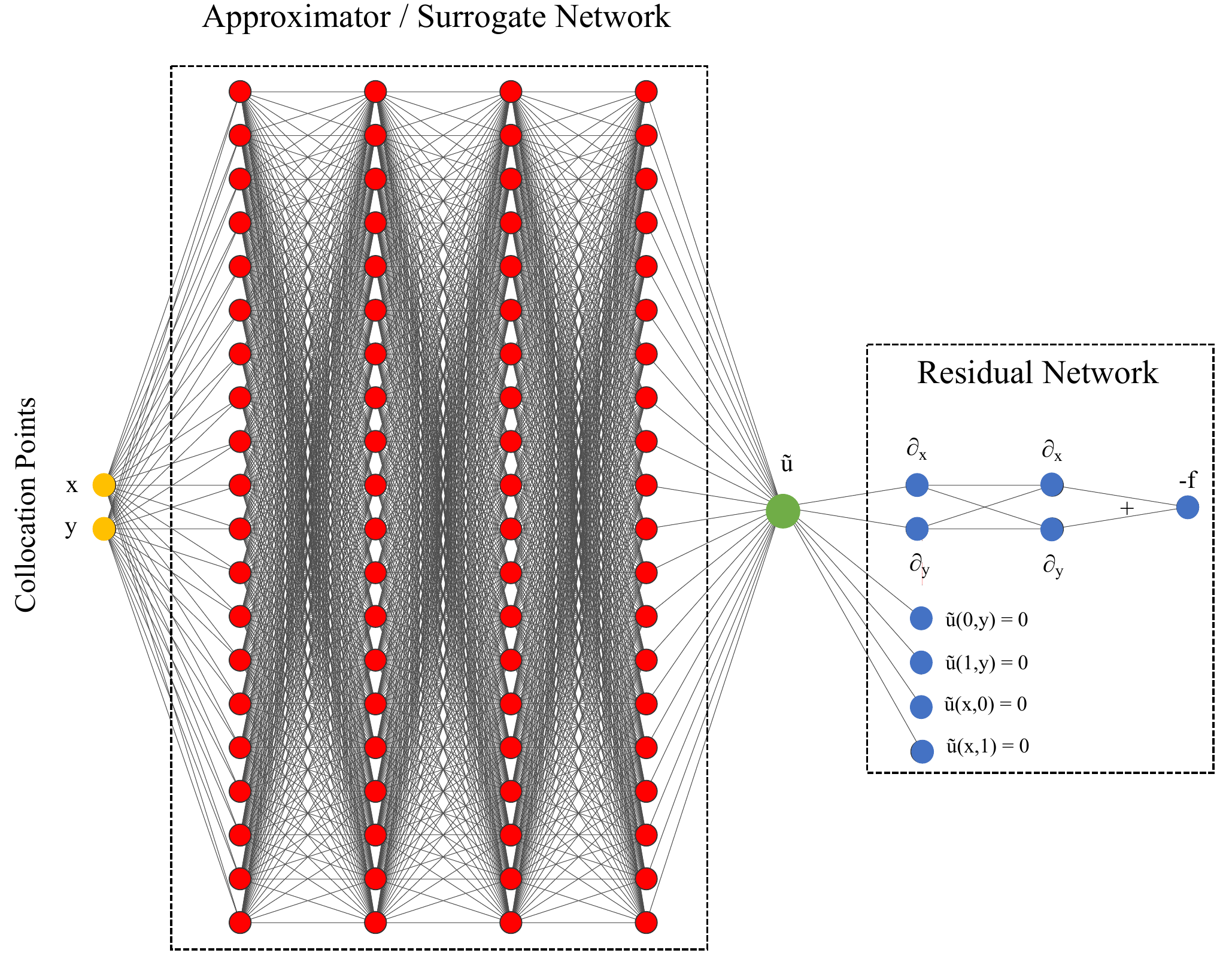}
\end{center}
\caption{A PINN to solve a Poisson problem $\partial^2_x u(x,y) + \partial^2_y u(x,y)= f(x,y)$ with associated Dirichlet boundary conditions. PINN consists of two basic interconnected networks. The first network (red vertices) provides a surrogate or approximation of the problem solution $u$. The network takes as input a point in the problem domain $(x,y)$ and provides an approximate solution $\tilde{u}$. This network weights and biases are trainable. The second network (blue vertices) takes the approximate solution from the first network and calculates the residual that is used as loss function to train the first network. The residual network includes the governing equations, boundary conditions and initial conditions (not included in the plot as the Poisson problem does not require initial conditions).}
\label{basicPINN}
\end{figure}

In its basic formulation, PINNs combine two networks together: an \emph{approximator} or \emph{surrogate} network and a residual network (see Figure~\ref{basicPINN})~\cite{raissi2019physics}. The approximator/surrogate network undergoes training and after it provides a solution $\tilde{u}$ at a given input point $(x,y)$, called \emph{collocation point}, in the simulation domain. The residual network encodes the governing equations and it is the distinctive feature of PINNs. The residual network is not trained and its only function is to provide the approximator/surrogate network with the residual (\emph{loss} function in DL terminology):
\begin{equation}
r =  \nabla^2 \tilde{u}(x,y) - f(x,y).
\label{resdi1}
\end{equation}
Differently from traditional methods often relying on finite difference approximation, the derivatives on the residual network graph, e.g, $\nabla^2 \tilde{u}(x,y)$ in Equation~\ref{resdi1}, are calculated using the so-called \emph{automatic differentiation}, or \texttt{autodiff}, that leverages the chain rule~\cite{baydin2018automatic} applied to the operations defined on the network nodes. In the solution of the Poisson Equation, the Laplacian operator is expressed as two successive first-oder derivatives of $\tilde{u}$ in the $x$ and $y$ directions and their summation (see the blue network nodes in Figure~\ref{basicPINN}).

In the inference/prediction phase, only the surrogate network is used to calculate the solution to the problem (remember that the residual network is only used in the training process to calculate the residual).

The approximator/surrogate network is a feedforward neural network~\cite{goodfellow2016deep}: it processes an input $x$ via $l$ layer of units (called also \emph{neurons}). The approximator/surrogate network expresses affine-linear maps ($Z$) between units and scalar non-linear activation functions ($a$) within the units:
\begin{equation}
\tilde{u} (x) = Z_l \circ a \circ Z_{l-1} \circ a ... \circ a \circ Z_{2} \circ a \circ Z_{1}(x).
\end{equation}
In DL, the most used activation functions are Rectified Linear Unit (\texttt{ReLU}), \texttt{tanh}, \texttt{swish}, \texttt{sine}, and \texttt{sigmoid} functions. See Ref.~\cite{ramachandran2017searching} for an overview of the different activation functions. As shown by Ref. ~\cite{mishra1}, PINNs requires sufficiently smooth activation functions. \textbf{PINNs with \texttt{ReLU} and other non-smooth activation functions, such as \texttt{ELU} and \texttt{SELU}~(Exponential and Scaled Exponential Linear Units) are not ``consistent/convergent" methods}: in the limit of an infinite training dataset a well-trained PINN with \texttt{ReLU}-like activation functions, the solution does not converge to the exact solution~\cite{mishra1}. This theoretical result is also confirmed by our experiments using \texttt{ReLU}-like activation functions. For this reason, we do not use \texttt{ReLU}-like activation functions in PINNs. 

The affine maps $Z$ are characterized by the weights and biases of the approximator/surrogate network:
\begin{equation}
Z_l x_l = W_l x_l + b_l,
\end{equation}
where $W_l$ is a \emph{weight} matrix for the layer $l$ and $b$ is the \emph{bias} vector. In PINNs, the weight values are initialized using the \emph{Xavier} (also called \emph{Glorot} when using the last name of the inventor instead) procedure~\cite{kumar2017weight}.

Typically, the PINN approximator/surrogate networks are fully connected networks consisting of 4-6 hidden layers(H) and 50-100 units per layer, similarly to the network in Figure~\ref{basicPINN}.  There are also successful experiments using convolutional and recurrent layers~\cite{gao2020phygeonet, nascimento2019fleet} but the vast majority of existing PINNs rely on fully-connected layers. In this work, we focus on studying the performance of fully-connected PINN.

The residual network is responsible for encoding the equation to solve and provide the loss function to the approximator network for the optimization process. In PINNs, we minimize the Mean Squared Error (MSE) of the residual (Equation~\ref{resdi1}):
\begin{equation}
MSE_r =  \frac {1}{N_{x_i,y_i}}\sum | r(x_i,y_i)  | ^2,
\label{trainingerror}
\end{equation}
where $N_{x_i,y_i}$ is the number of collocation points. \textbf{In PINNs, the collocation points constitute the training dataset}. Note that $MSE_r$ depends on the size of the training of the dataset ($N_{x_i,y_i}$), e.g., the number of collocation points. In practice, a larger number of collocation points leads to an increased MSE value. $MSE_r$ depends also on on the distribution of our collocation points. The three most used dataset distributions are: \texttt{uniform} (the dataset is uniformly spaced on the simulation domain as on a uniform grid), \texttt{pseudo-random} (collocations points are sampled using pseudo-random number generator) and \texttt{Sobol} (collocation points are from the Sobol low-discrepancy sequence). Typically, the default training distribution for PINNs is \texttt{Sobol}, like in quasi-Montecarlo methods.

Recently, several PINN architectures have been proposed. PINNs differentiate on how the residual network is defined.  For instance, \texttt{fPINN} (fractional PINN) is a PINN with a residual network capable of calculating residuals of governing equations including fractional calculus operators \cite{pang2019fpinns}.  \texttt{fPINN}  combines automatic differentiation with numerical discretization for the fractional operators in the residual network. \texttt{fPINN} extends PINN to solve integral and differential-integral equations. Another important PINN is \texttt{vPINN}  (variational PINN): they include a residual network that uses the variational form of the problem into the loss function~\cite{kharazmi2019variational} and an additional shallow network using trial functions and polynomials and trigonometric functions as test functions. A major advantage with respect to basic PINNs is that in the analytical calculation by integrating by parts the integrand in the variational form, we can the order of the differential operators represented by the neural networks, speeding up the training and increasing PINN accuracy. \texttt{hp-VPINN} is an extension of \texttt{vPINN} that allows hp-refinement via domain decomposition as h-refinement and projection onto space of high order polynomials as p-refinement~\cite{kharazmi2020hp}. In this work, we use the original residual network as shown in Figure~\ref{basicPINN}.

In the training phase, an optimization process targeting the residual minimization determines the weights and biases of the surrogate network. Typically, we use two optimizers in succession: the Adam optimizer as first and then a Broyden-Fletcher-Goldfarb-Shanno (BFGS) optimizer~\cite{fletcher2013practical}. BFGS uses the Hessian matrix (curvature in highly dimensional space) to calculate the optimization direction and provides more accurate results. However, if used directly without using the Adam optimizer can rapidly converge to a local minimum (for the residual) without exiting. For this reason, the Adam optimizer is used first to avoid local minima, and then the solution is refined by BFGS. We note that the typical BFGS used in PINNs is the L-BFGS-B: L-BFGS is a limited-memory version of BFGS to handle problems with many variables, such as DL problems; the BFGS-B is a variant of BFGS for bound constrained optimization problems. In our work, we tested several optimizers, including Newton and Powell methods, and found that L-BFGS-B provides by far the highest accuracy and faster convergence in all our test problems. \textbf{L-BFGS-B is currently the most critical technology for PINNs}.

An \emph{epoch} comprises all the optimizer iterations to cover all the datasets. In PINNs, typically, thousands of epochs are required to achieve accurate results. By nature, PINNs are under-fitted: the network is not complex enough to accurately capture relationships between the collocation points and solution. Therefore, an extensive dataset increase improves the PINN performance; however, the computational cost increases raising the data set size. 

One crucial point related to PINNs is whether a neural network can approximate simultaneously and uniformly the solution function and its partial derivatives. Ref.~\cite{lu2019deepxde} shows that feed-forward neural nets with enough neurons can achieve this task. A formal analysis of the errors in PINNs is presented in Refs. \cite{mishra1,lu2019deepxde}.  

\textbf{An important fact determining the convergence behavior of the DL networks and PINN linear solvers is the Frequency-principle (F-principle)}:  \emph{DNNs often fit target functions from low to high frequencies during the training process} \cite{xu2019frequency}. The F-principle implies that in PINNs, the low frequency / large scale features of the solution emerge first, while it will take several training epochs to recover high frequency / small-scale features. This 

Despite the recent introduction of PINNs, several PINN frameworks for PDE solutions exist. All the major PINN frameworks are written in Python and rely either on \texttt{TensorFlow}~\cite{abadi2016tensorflow} or \texttt{PyTorch}~\cite{paszke2019pytorch} to express the neural network architecture and exploit auto-differentiation used in the residual network. Together with \texttt{TensorFlow}, \texttt{SciPy}~\cite{virtanen2020scipy}  is often used to use high-order optimizers such as L-BFGS-B. Two valuable PINN Domain-Specific Languages (DSL) are \texttt{DeepXDE}~\cite{lu2019deepxde} and \texttt{sciANN}~\cite{haghighat2020sciann}. DeepXDE is an highly customizable framework with TensorFlow 1 and 2 backend and it supports basic and fractional PINNs in complex geometries. \texttt{sciANN} is a DSL based on and similar to \texttt{Keras}~\cite{gulli2017deep}. In this work, we use the \texttt{DeepXDE} DSL.

\subsection{An Example: Solving the 2D Poisson Equation with PINN}\label{sec:poisson}
To showcase how PINNs work and provide a baseline performance in terms of accuracy and computational cost, we solve a Poisson problem in the unit square domain with a source term $f(x,y)$ that is smooth, e.g., differentiable, and contains four increasing frequencies:
\begin{equation}
f(x,y) = \frac{1}{4}  \sum_{k=1}^4 (-1)^{k+1} 2 k \sin(k  \pi x) \sin( k  \pi y).
\label{manysin} 
\end{equation}
We choose such a source term as it has a simple solution and to show the F-principle's impact on the convergence of PINN to the numerical solution: we expect the lower frequency components, e.g., $k=1$, to convergence faster than the higher frequency components present in the solution ($k = 2, 3, 4$).

We use a fully-connected four-layer PINN with a \texttt{tanh} activation function for the approximator/surrogate network for demonstration purposes and without a loss of generality. The input layer consists of two neurons (the $x$ and $y$ coordinates of one collocation point), while each hidden and output layers comprise 50 neurons and one neuron, respectively. The weights of the network are initialized with the Xavier method. As a reminder, the approximator/surrogate network's output is the approximate solution to our problem. The residual network is a graph encoding the Poisson equation and source term and provides the loss function (Equation~\ref{trainingerror}) to drive the approximator/surrogate network's optimization. At each, a collocation point within the problem domain is drawn from the \texttt{Sobol} sequence. The training data set consists of 128 $\times$128 collocation points on the domain and additional 4,000 collocation points on the boundary for a total of 20,384 points. We train the approximator/surrogate network 10,000 of Adam optimizer epochs with a learning rate $\lambda$ equal to 0.001 (the magnitude of the optimizer vector along the direction to minimize the residual), followed by 13,000 epochs of L-BFGS-B optimizer. We use the \texttt{DeepXDE} DSL for our PINN implementation.

\begin{figure}[h!]
\begin{center}
\includegraphics[width=\textwidth]{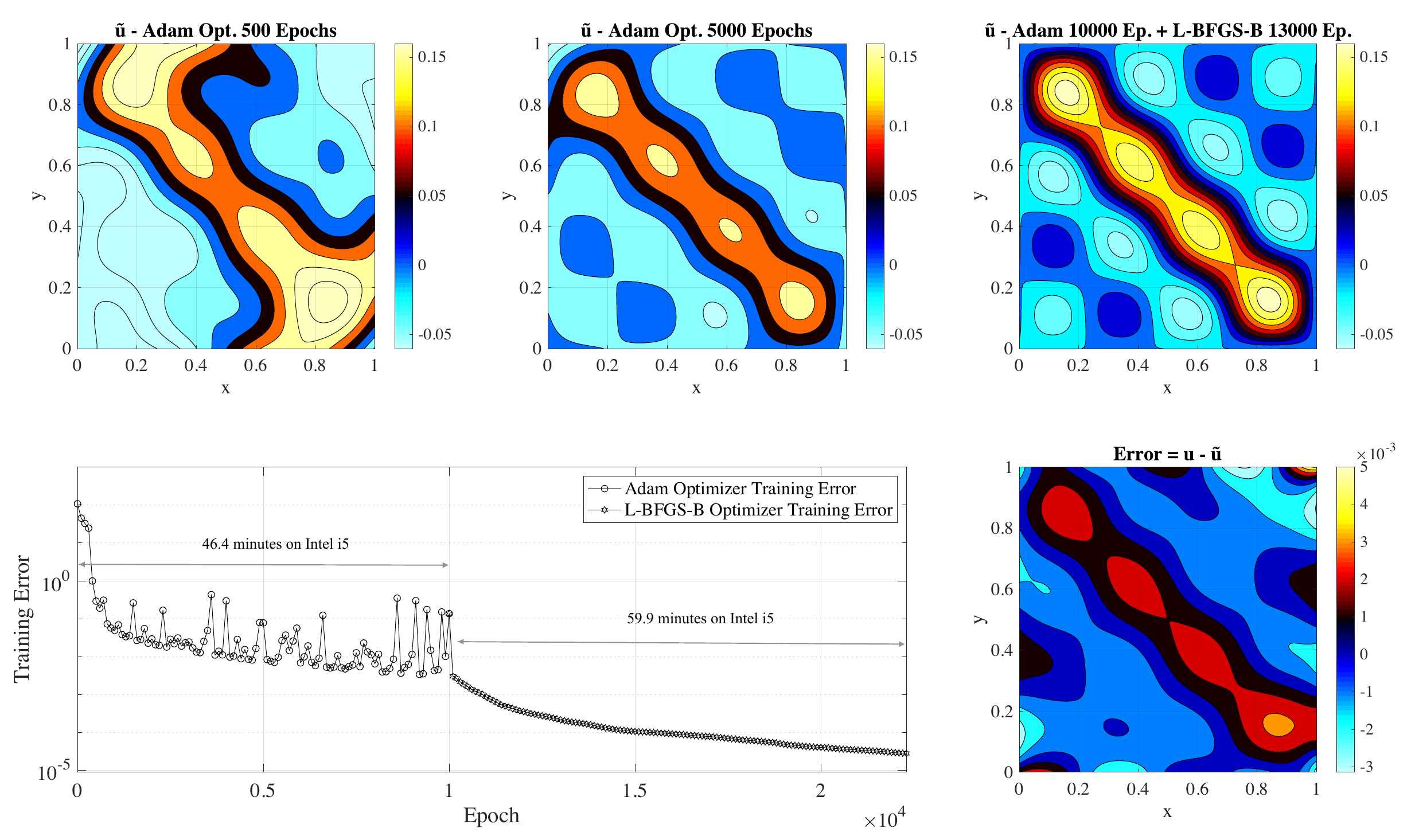}
\end{center}
\caption{The top panels show the solution of the Poisson equation at different epochs using a PINN. The bottom panel shows the training error for an initial training with Adam's optimizer (10,000 epochs), followed by L-BFGS-B (13,000 epochs). The plot also includes the total time for training the PINN on a dual-core Intel i5 processor. The right bottom subplot presents the error of the final solution compared to the exact solution.}
\label{basicPINNresults}
\end{figure}

Figure \ref{basicPINNresults} shows the Poisson equation's approximate solution with the source term of Equation~\ref{manysin} at different epochs, the training error, and the error of the PINN solution after the training is completed. The Figure \ref{basicPINNresults} top panels present the contour plot of the approximator/surrogate solution on a 128$\times$128 uniform grid after 500, 5,000 and 23,000 epochs. To determine the solution at each epoch, we take the approximate/surrogate network and perform inference/prediction using the points of the 128$\times$128 uniform grid. By analyzing the approximate solutions' evolution (top panels of Figure~\ref{basicPINNresults}), it is clear that the PINN resolves the low-frequency component present in the solution: a yellow band appears along the diagonal of the plot while local peaks (small islands in the contour plot) are not resolved. As the training progresses, localized peaks associated with the source term's high-frequencies appear and are resolved. The bottom right panel of Figure~\ref{basicPINNresults} shows a contour plot of the error after the training is completed. The maximum pointwise error is approximately 5E-3. We note that a large part of the error is located in the proximity of the boundaries. This issue results from the \emph{vanishing-gradient} problem \cite{wang2020understanding}: unbalanced gradients back-propagate during the model training. This issue is similar to the numerical \emph{stiffness} problem when using traditional numerical approaches. One of the effective technique to mitigate the \emph{vanishing-gradient} problem is to employ locally (to the layers or the node) adaptive activation functions~\cite{jagtap2020locally}. Additional techniques for mitigating \emph{vanishing-gradient} problem are the usage of ReLU activations functions and batch normalization.

The bottom panel of Figure \ref{basicPINNresults} shows the training error's evolution calculated with Equation~\ref{trainingerror}. In this case, the initial error is approximately 1.08E2 and decreases up to 2.79E-5 at the end of the training. The initial error mainly depends on the training data set size: small input data sets reduce training error that does not translate to higher accuracy in the solution of the problem. However, the training is a reasonable metric when comparing the PINN performance when using the same data set size.

By analyzing the evolution of the training error, it is clear that the Adam optimizer training error stabilizes approximately in the range of 5E-3 - 1E-2 after 2,000 epochs, and we do not observe any evident improvement after 2,000 epochs of Adam optimization. The L-BFGS-B optimizer leads the error from 5E-3 - 1E-2 to 2.79E-5 and is responsible for the major decrease of the training error. However, we remind that L-BFGS-B is not used at the beginning of the training as it can converge quickly to a wrong solution (a local minimum in the optimization problem).

To provide an idea of the PINN training's overall computation cost, we also report the total time for training the PINN in this basic non-optimized configuration on a dual-core Intel i5 2.9 GHz CPU. The total training execution time is 6,380 seconds, corresponding to approximately 1.5 hours. For comparison, the solution of the same problem with a uniform grid size 128$\times$128 on the same system with the \texttt{petsc4py} CG solver~\cite{dalcin2011parallel,balay2019petsc} requires 92.28 seconds to converge to double-precision machine epsilon. Basic PINN's direct usage to solve the Poisson problem is limited for scientific application given the computational cost and the relatively low accuracy. In the next sections, we investigate which factors impact the PINN performance and its accuracy. We design a PINN-based solver to have comparable performance to state-of-the-art linear solvers such as \texttt{petsc4py}.

\section{Characterizing PINNs as Linear Solvers}\label{sec:tune}
To characterize the PINNs performance for solving the Poisson equation, we perform several parametric studies varying the approximator/surrogate network size, activation functions, and training data size and distribution. We also investigate the performance enhancement achieved by using the transfer learning technique to initialize with the network weights obtained solving the Poisson equation with a different source term~\cite{weiss2016survey}. During our experiments, we found that two relatively different configurations of the network are required in the case of the source term of the Poisson equation is smooth on non smooth, e.g. non-differentiable. For this reason, we choose two main use cases to showcase the impact of different parameters. For the smooth source term case, we take the source term from Equation~\ref{manysin} (the example we showcased in the previous section). For the non-smooth source term case, we take a source term that is zero everywhere except for the points enclosed in the circle, centered in $(0.5,0.5)$ with radius $0.2$: 
\begin{equation}
f(x,y) = 1 \; \textnormal{for} \; \sqrt{(x-0.5)^2  + (y-0.5)^2 } \leq 0.2. 
\label{nonsmooth}
\end{equation}
As baseline configuration, we adopt the same configuration described in the previous section: a fully-connected network with four hidden layers of 50 units, and \texttt{tanh} activation function. The data set consists of 128$\times$128 collocation points in the domain and 4,000 points on the boundary. Differently from the previous configuration, we reduce the training epochs to 2,000 for the Adam Optimizer (the training error do not decrease after 2,000 epochs) and 5,000 for the L-BFGS-B optimizer.

\begin{figure}[h!]
\begin{center}
\includegraphics[width=\textwidth]{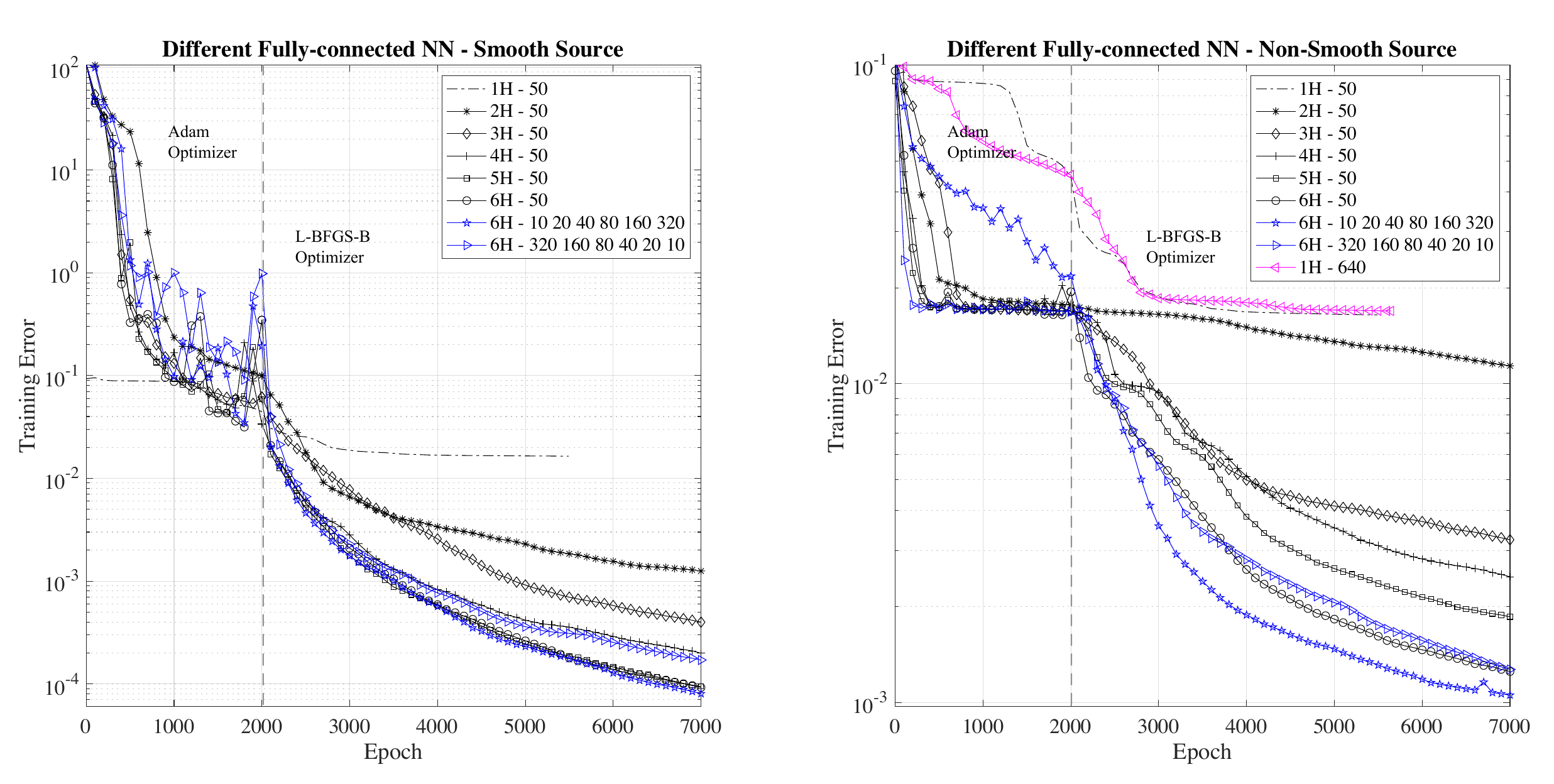}
\end{center}
\caption{Training error for different fully-connected PINN depth: one (\texttt{1H}), two (\texttt{2H}), three (\texttt{3H}), four (\texttt{4H}), five (\texttt{5H}) and six (\texttt{6H}) hidden layers with 50 neurons each. We also consider the training error for PINNs with six hidden layers and \texttt{10-20-40-80-160-320} and \texttt{320-160-80-40-20-10} units per hidden layer, respectively.}
\label{architecture}
\end{figure}

The first experiments we perform is to evaluate the impact of the network size (depth and units per layer) on the training error. To understand the impact of surrogate neural network depth, we perform training with layers of 50 neurons with one (\texttt{1H}), two (\texttt{2H}), three (\texttt{3H}), four (\texttt{4H}), five (\texttt{5H}) and six (\texttt{6H}) hidden layers (\texttt{H} stands for hidden layer). We present the evolution of training error in Figure \ref{architecture}. By analyzing this figure, it is clear that shallow networks consisting of one or two hidden layers do not perform, and the PINN learning is bound in learning after few thousand epochs. Even one layer with large number of units, e.g., one hidden layer with 640 units (see the magenta line in the right panel of Figure \ref{architecture}), do not lead to better performance as demonstration that depth is more important than breadth in PINN. Deeper networks with more than three layers lead to lower final training errors and improved learning. However, we find that the final training error saturates for PINNs with more than six hidden layers (results not shown here) for the two test cases. An important aspect for the deployment of PINN in scientific applications is that the performance of PINNs with four and more hidden layers have comparable performance in the first 500 epochs of the Adam and L-BFGS-B optimizers. Taking in account that the PINN computational cost for PINNs increases with the number layers and realistically only few hundred epochs are necessary for PINN to be competitive with HPC solvers, PINNs with four hidden layers provide the best trade-off in terms of accuracy and computational performance. 

For the six hidden layers case, we also check the importance of having a large/small number of units at the beginning/end of the network: we consider the performance of PINN with six hidden layers and \texttt{10-20-40-80-160-320} and \texttt{320-160-80-40-20-10} units per hidden layer, respectively. We find that to have a large number of units at the beginning of the network and small number of units at the end of the network is detrimental to the PINN performance (a six hidden layer network in this configuration has the same performance of a five hidden layer PINN). Instead, to have a small number of units at the beginning of the network and a large number of units at the end of the network is beneficial to the PINN. \textbf{This observation hints that initial hidden layers might responsible for encoding the low-frequencies components (fewer points are needed to represent low-frequency signals) and the following hidden layers are responsible for representing higher-frequency components (several points are needed to represent high-frequency signals)}. However, more experiments are needed to confirm this hypothesis.

\textbf{The most impactful parameter for achieving a low training error is the activation function}. This fact is expected as activation functions are nothing else than non-linear interpolation functions (similarly to nodal functions in FEM): some interpolation function might be a better fit to represent the different source terms. For instance, sigmoid functions are a good fit to represent non-differentiable source terms exhibiting discontinuities. On the contrary, a smooth \texttt{tanh} activation function can closely represent smooth functions.

\begin{figure}[h!]
\begin{center}
\includegraphics[width=\textwidth]{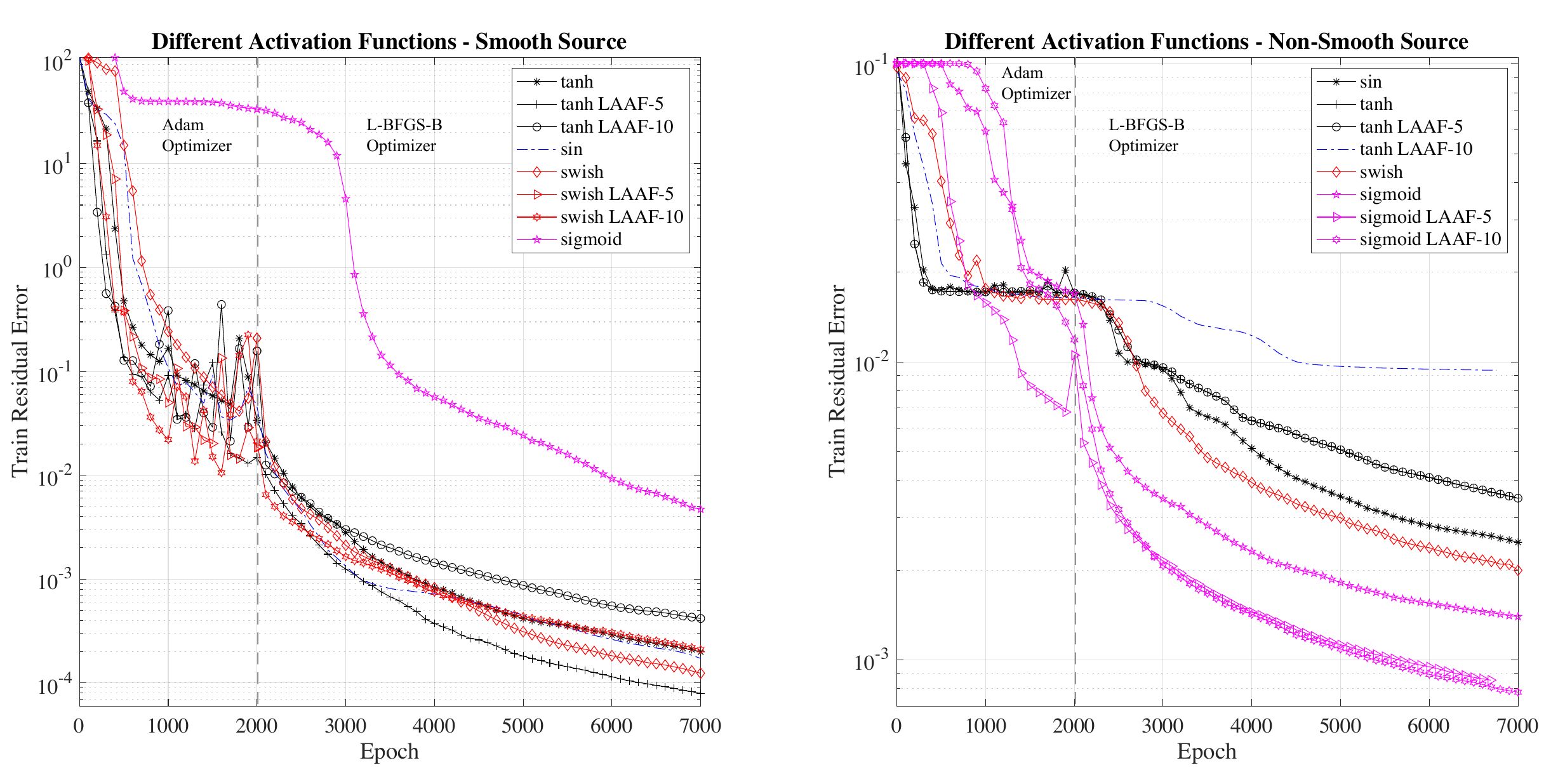}
\end{center}
\caption{Training error for different activation functions. The two test cases show rather different performance: the best activation function for smooth source term case is \texttt{tanh}, while it is \texttt{sigmoid} for the non-smooth source term case. Local (to the layer) adaptive activation functions provide a reduction of the training error.}
\label{activation}
\end{figure}
We investigate the impact of different activation functions and show the evolution of the training errors in Figure~\ref{activation}. Together with traditional activation function, we also consider the Locally Adaptive Activation Functions (\texttt{LAAF}): with this technique, a scalable parameter is introduced in each layer separately, and then optimized with a variant of stochastic gradient descent algorithm~\cite{jagtap2020locally}. The LAAF are provided in the \texttt{DeepXDE} DSL. We investigate LAAF with factor of 5 (\texttt{LAAF-5}) and 10 (\texttt{LAAF-10}) for the \texttt{tanh}, \texttt{swish} and \texttt{sigmoid} cases. The \texttt{LAAF} usage is critical to mitigate the \emph{vanishing-gradient} problem.

The activation function's different impact for the two test cases (smooth and non-smooth source terms) is clear when analyzing the results presented in Figure~\ref{activation}. In the smooth source term case, the best activation function is the locally (to the layer) adaptive \texttt{tanh} activation function with factor 5 (\texttt{LAAF5 - tanh}). In the case of the non-smooth source term, the \texttt{sigmoid} activation function outperforms all the other activation functions. In particular, in this case, the best activation function is the locally (to the layer) adaptive sigmoid activation function with factor 10 (\texttt{LAAF10 - sigmoid}).

\begin{figure}[h!]
\begin{center}
\includegraphics[width=\textwidth]{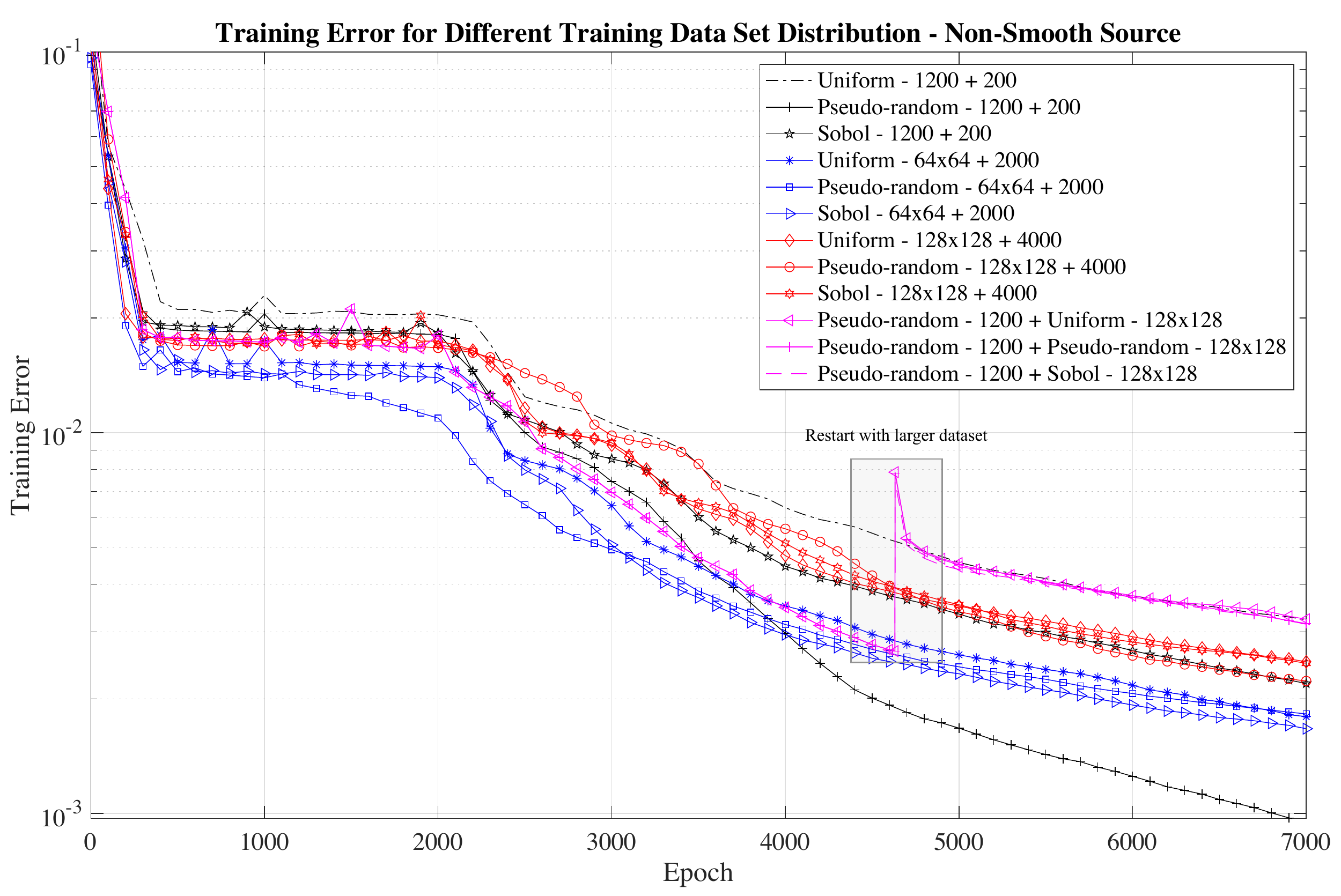}
\end{center}
\caption{Training error for different data set (1,200 points in the domain and 200 on the boundary, 64$\times$64 points in the domain and 2,000 on the boundary, 128$\times$128 points in the domain and 4,000 on the boundary) and different distribution (\texttt{uniform}, \texttt{pseudo-random} and \texttt{Sobol}).}
\label{impactdata}
\end{figure}

As we mentioned in Section~\ref{sec:poisson}, the data size impacts the training errors. Large data sets increase the PINN accuracy but have larger training errors than the training with small data set because of the error definition (see Equation \ref{trainingerror}). For this reason, the training error should be compared only for training using the same training data set size.  We investigate the impact of three different input data size ( 1- 1,200 points in the domain and 200 on the boundary, 2-64$\times$64 points in the domain and 2,000 on the boundary, 3- 128$\times$128 points in the domain and 4,000 on the boundary) with three collocation point distributions (\texttt{uniform}, \texttt{pseudo-random}, and \texttt{Sobol} sequence) for the non-smooth source term. We show the results in Figure~\ref{impactdata}. 

In general, we find that the collocation point distribution does not have a considerable impact on the training error for large data sets: the \texttt{Sobol} and \texttt{pseudo-random} distributions have a slightly better performance than the \texttt{uniform} distribution. For small data sets, \texttt{pseudo-random} distribution result in lower training errors. 

We also study the impact of having a \emph{restart} procedure: we train first the PINN with a small data set 1,200 points in the domain and 200 on the boundary) for 4,500 epochs (and then re-train the same network with a large data set (128$\times$128 points in the domain and 4,000 on the boundary) for 2,500 cycles (see the magenta lines and the grey box in Figure~\ref{impactdata}). Such a restart capability would lead to a large computational saving. However, the results show that to retrain with a large data set does not lead to a decreased error and result in the highest training error.

\section{The Importance of Transfer Learning}\label{sec:transferlearning}

In this study, we found that the usage transfer learning technique is critical for training PINNs with a reduced number of epochs and computational cost. The transfer learning technique consists of training a network solving the Poisson equation with a different source term.  We can then initialize the PINN network we intend to solve with the first fully trained network weights and biases. In this way, the first PINN \emph{transfers} the learned information about encoding to the second PINN. To show the advantage of transfer learning in PINN, we solve two additional test cases with smooth and non-smooth source terms. For the test case with the smooth source term, we solve the Poisson equation with source term $f(x,y) = 10(x(x - 1) + y(y - 1)) -2\sin (\pi x)\sin (\pi y) + 5(2\pi x)\sin (2\pi y)$. 

\begin{figure}[h!]
\begin{center}
\includegraphics[width=\textwidth]{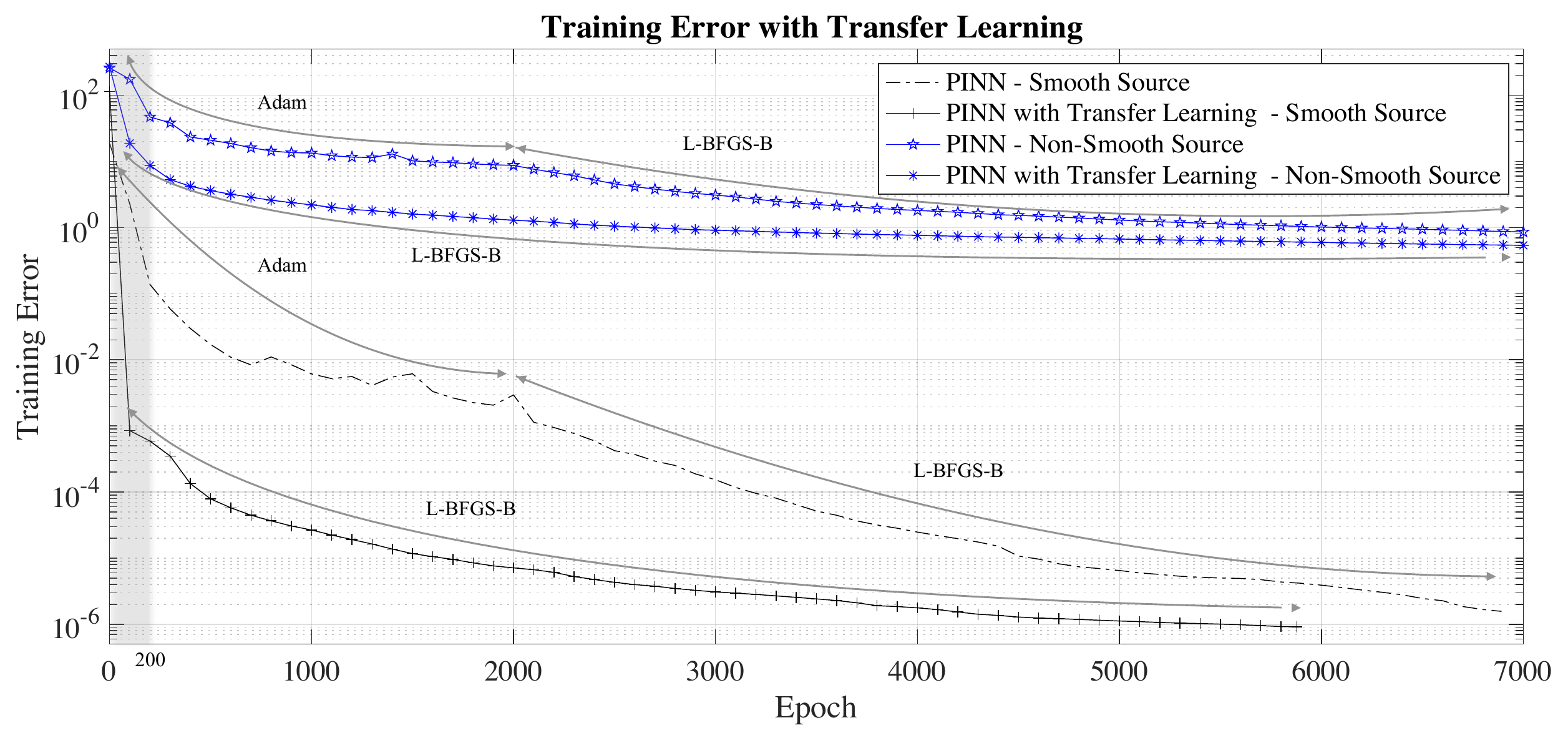}
\end{center}
\caption{Training error with and without transfer learning for the smooth and non-smooth source test cases.}
\label{transferlearning}
\end{figure}

We initialize the network with the results obtained during the training with Equation~\ref{manysin} as a source term. One of the major advantages of transfer-learning is that we can start the L-BFGS-B optimizer after very few Adam solvers epochs (empirically,we found that 10 Adam epochs ensure that L-BFGS-B will avoid local minima). L-BFGS-B has faster convergence than the Adam optimizer and therefore the training is quicker. When not using transfer-learning, we train the PINN with 2,000 epochs of Adam optimizer, followed by 5,000 epochs of L-BFGS-B. When using L-BFGS-B, we perform 10 epochs of Adam optimizer, followed by 6,955 L-BFGS-B epochs.

The black lines in Figure~\ref{transferlearning} show a comparison of the training error for a network initialized with Xavier weight initialization, e.g., without transfer learning ($-.$ black line) and with transfer learning ($-+$ black line). In this case, transfer learning usage allows gaining two orders of improvement in the training error in less than 1,000 epochs.

For the test case with non-smooth source term, we introduce and additional test case solving the Poisson equation with a source term that is everywhere zero except in a circle with radius $0.1$ and centered in the $x$ and $y$ coordinates (0.7,0.7).
\begin{equation}
f(x,y) = - 10 \; \textnormal{for} \;  \sqrt{(x-0.7)^2  + (y-0.7)^2 } \leq 0.1. 
\label{nonsmooth2}
\end{equation}
For transfer learning, we use the PINN weights obtained training the network to solve the Poisson equation with source term of Equation~\ref{nonsmooth2}. The blue lines in Figure~\ref{transferlearning} are the training error without transfer learning. As in the case of smooth-source term, the usage of transfer learning rapidly decreases the training error.

We note that usage of the transfer learning leads to an initial (less than 200 L-BFGS-B epochs) \emph{super-convergence} to a relatively low training error. For this reason, \textbf{transfer-learning is a necessary operation to make PINN competitive with other solvers used in scientific computing}.

The major challenge for using transfer-learning is to determine which pre-trained PINN to use. In simulation codes, solving the same equation with different source term at each time step, an obvious choice is a PINN that solves the governing equations with a source term at one of the time step. For other cases, we found that PINNs solving problems with source terms containing high-frequency components (possibly more than one component) are suitable for transfer-learning in general situations. We also found that PINNs solving problem with only one low frequency component as source term are not beneficial for transfer learning: their performance is equivalent to the case without transfer learning.

\section{The Old and the New: Integrating PINNs into Traditional Linear Solvers}\label{sec:integrate}
In Section \ref{sec:poisson}, we observed that direct usage of PINN to solve the Poisson equation is still limited by the large number of epochs required to achieve an acceptable precision. One possibility to improve the performance of PINN is to combine PINN with traditional iterative solvers such as the Jacobi, Gauss-Seidel and multigrid solvers~\cite{quarteroni2010numerical}.

PINN solvers' advantage is the quick convergence to the solution's low frequencies components. However, the convergence to high-frequency features is slow and requires an increasing number of training iteration/epochs. This fact is a result of the the F-principle. Because of this, PINNs are of limited usage when the application requires highly accurate solutions. As suggested by Ref.~\cite{xu2019frequency}, in such cases, the most viable option is to combine PINN solvers with traditional solvers that can converge rapidly to the solution's high-frequency components (but have low convergence for the low-frequency components). Such methods introduce a computational grid and we compute the differential operators with a finite difference scheme. In this work, we choose the Gauss-Seidel method as it exhibits higher convergence rate than the Jacobi method. Each  Gauss-Seidel solver iteration for solving the Poisson equation (Equation \ref{poisson}) is:
\begin{equation}
u_{i,j}^{n+1} = 1/4 (u_{i+1,j}^{n} + u_{i-1,j}^{n+1} + u_{i,j+1}^{n}  + u_{i,j-1}^{n+1}  - \Delta x  \Delta y f_{i,j}),
\label{GSeq}
\end{equation}
where $i$ and $j$ are the cell index, $\Delta x$ and $\Delta y$ are the grid cell size in the $x$ and $y$ direction, and $n$ is the iteration number. Usually, the Gauss-Seidel method stops iterating when $||u_{n+1} - u^n ||_2 \leq \delta$, where $|| ... ||$ is the Euclidean norm and $ \delta$ is a so-called tolerance and it is chosen as an arbitrarily small value.

Both the Jacobi and Gauss-Seidel methods show fast convergence for small-scale features: this is because the update of unknown values involves only the values of the neighbor points (stencil defined by the discretization of a differential operator). Between two different iterations, the information can only propagate to neighbor cells.

\begin{figure}[bt]
\begin{center}
\includegraphics[width=0.4\textwidth]{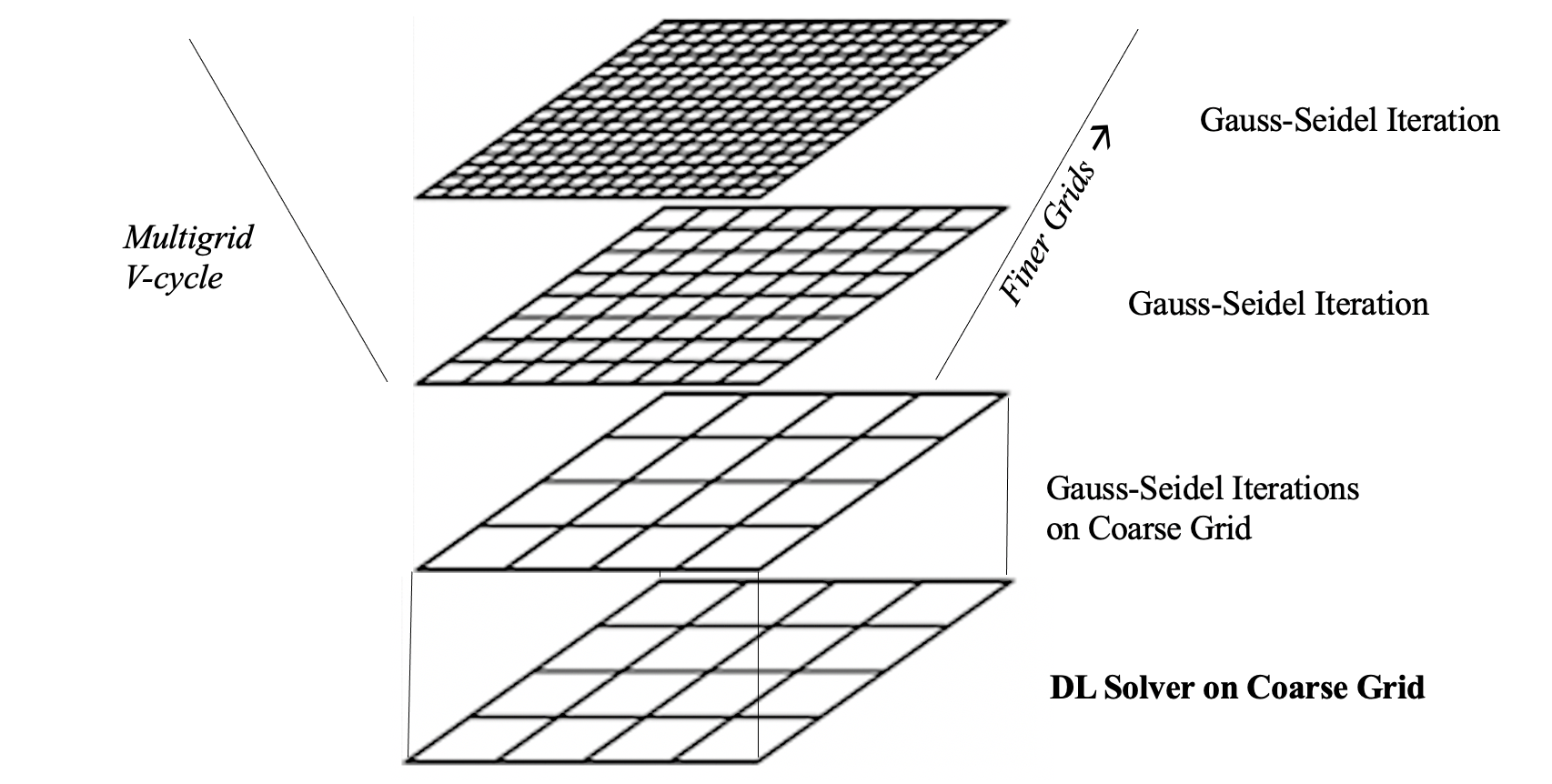}
\end{center}
\caption{The hybrid solvers relies on the DL linear solver to determine the solution on a coarse grid that is refined through a multigrid V-cycle performing Gauss-Seidel iterations on finer grids.}
\label{vcycle}
\end{figure}

In this work, we combine traditional approaches with new emerging DL methods as shown in Figure~\ref{vcycle}. Overall, the new solver consists of three phases. We use first the DL PINN solver to calculate the solution on a coarse grid. As second phase, we refine the solution with Gauss-Seidel iterations on the coarse grid until a stopping criteria is satisfied. The third phase is a multigrid V-cycle: we linearly interpolate (or \emph{prolongate} in multigrid terminology) to finer grids and perform a Gauss-Seidel iteration for each finer grid. In fact, several multigrid strategies with different level of sophistications can be sought. However, in this work we focus on a very simple multigrid approach, based on the Gauss-Seidel method and linear interpolation across different grids. The crucial point is that we train a PINN to calculate the solution of the problem on the coarse grid, replacing the multigrid \emph{restriction} (or \emph{injection}) steps in just one phase.

\begin{figure}[h!]
\begin{center}
\includegraphics[width=0.85\textwidth]{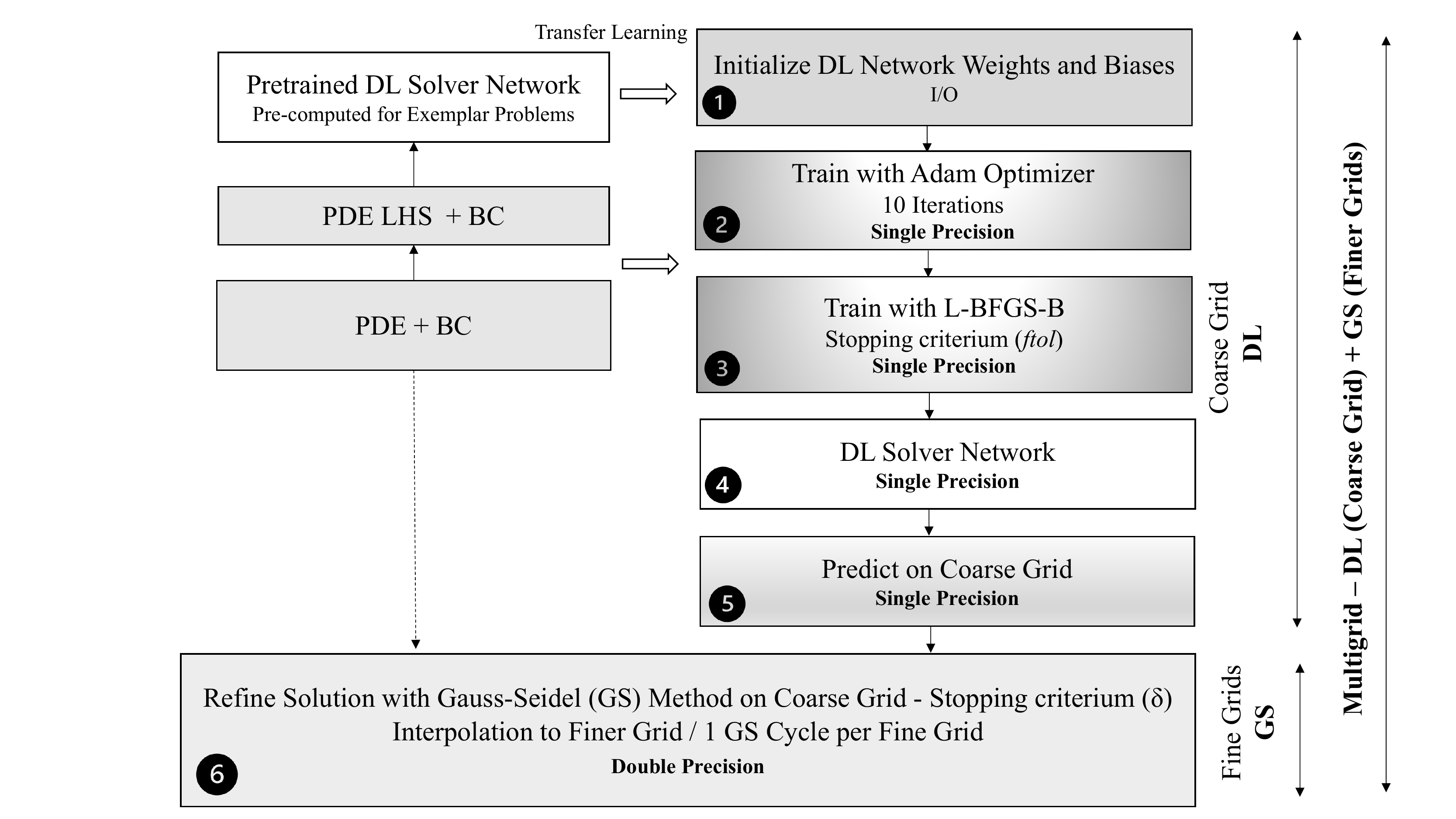}
\end{center}
\caption{Structure of the hybrid multigrid solver combining the DL and Gauss-Seidel solvers. Pre-trained networks are pre-computed and used to initialize the DL network. Two main parameters $ftol, \delta$  determine the accuracy and the performance of the hybrid solver.}
\label{hybridpinn}
\end{figure}

Figure~\ref{hybridpinn} shows a more detailed diagram of a hybrid multigrid solver combining a DL solver to calculate the solution on a coarse grid with a Gauss-Seidel solver to refine the solution and interpolate to finer grid. Because the DL solver convergences quickly to the low-frequency  coarse-grained components of the solution while high-frequency small-scale components of the solutions are not accurately solved, we perform the training in single-precision floating-point. This would speed-up the training on GPUs (not used in this work) where the number of single-precision floating-point units (FPUs) is higher than CPU.

The hybrid DL solver comprises six basic steps, represented in Figure~\ref{hybridpinn} :
\begin{enumerate}
\item \textbf{Initialize the network weights and biases} - We load from the disk the network structure and initialize the network. To accelerate the convergence, we rely on transfer-learning: we train a network to solve a similar problem and initialize the network. It is important that the same governing equations, boundary conditions and architecture are used. The weights and biases are in single floating-point precision. The time for completing this step is negligible with respect to the total time of the hybrid solver.
\item \textbf{Train with Adam Optimizer (10 Epochs)} - We run the Adam optimizer just for a short number of epochs to avoid the consequent L-BFGS-B optimizer converging quickly to the wrong solution (local minimum). By running several tests, we found empirically that only 10 Adams epochs are needed to avoid L-BFGS-B optimizer to converge to the wrong solution. The time for completing this step is typically negligible.
\item \textbf{Train with L-BFGS-B Optimizer} - We run the training with the L-BFGS-B optimizer. The stopping criterium is determined by the \emph{ftol} parameter: the training stops when $(r_k - r_{k+1})/\max(|r_k|,|r_{k+1}|,1) \leq  ftol$, where $k$ is the iteration of the optimizer and $r$ is the value of the function to be optimized (in our case the residual function). Typically, the time for completing the L-BFGS-B dominates is a large part of the execution time of the hybrid solver. To compete with traditional approaches for solving Poisson equation, we set a maximum number of epochs to 1,000. 
\item  \textbf{DL solver is obtained at the end of the training process} - The solver can inference the solution at given collocation points or save it for future transfer-learning tasks, e.g., a simulation repeats the computation of the Poisson equation at different time steps.
\item  \textbf{The Approximator/Surrogate Network is used to calculate the solution on the coarse grid of the multigrid solver} - We calculate the solution of our problem on the coarse grid of a multigrid solver. This operation is carried with single-precision floating point numbers since high-accuracy is not needed in this step. The result is then cast to double precision for the successive Gauss-Seidel solver. This inference computational time is typically negligible when compared to the total execution time.
\item \textbf{Refine the solution with the Gauss-Seidel Method on the coarse grid and interpolate on fine grids} - We perform first Gauss-Seidel iterations to refine the solution on the coarse grid. This solution refinement is critical to remove the vanishing-gradient problem at the boundary. The Gauss-Seidel iteration on the coarse grid stops when $||u^{n+1} - u^n ||_2 \leq \delta$ where $n$ is the iteration number. After the Gauss-Seidel method stops on the coarse grid, linear interpolation to finer grids and a Gauss-Seidel iteration per grid are computed. As example, to solve the problem on a 512$\times$512 grid, we perform the following steps: 1) use the DL solver to calculate the solution on 64$\times$64 grid; 2) refine the solution with the Gauss-Seidel method on the 64$\times$64 grid until convergence is reached; 3) carry out a linear interpolation to the 128x128 grid;  4) perform a Gauss-Seidel iteration on the 128$\times$128 grid; 5) carry out a linear interpolation to 256$\times$256 grid; 6) perform a Gauss-Seidel iteration on the 256$\times$256 grid; 7) carry out a linear interpolation to 512$\times$512 grid; 8) perform a final Gauss-Seidel iteration on the 512$\times$512 grid. The interpolation and Gauss-Seidel iterations corresponds to the V-cycle in the multigrid method as shown in Figure~\ref{vcycle}.
\end{enumerate}

We test the hybrid modified solver against the same problem shown in Section \ref{sec:poisson}: we solve the Poisson equation with source term of Equation~\ref{manysin}. Leveraging the knowledge gained in the characterization study of Section~\ref{sec:tune}, we use a four hidden layer fully-connected neural network with 50 neurons per hidden layer. To optimize the convergence for solving the Poisson equation with a smooth source term, we rely on \texttt{LAAF-5 tanh} activation functions: these activations functions provided the best performance in our characterization study. For the transfer learning, we pre-train a network for 2,000 Adam optimizer epochs and 5,000  L-BFGS-B optimizer epochs to solve a Poisson equation with a source term equal to $-2 \sin( \pi x) \sin(\pi  y) - 72 \sin(6 \pi x ) \sin (6 \pi y)$. We use an input data set consisting of 100$\times$100 points in the integration domain and 2,000 points on the boundaries for the DL solver. We use the \texttt{Sobol} sequence as training data set distribution. The network weights and biases for transfer learning are saved as checkpoint / restart files in TensorFlow.

\begin{figure}[bt]
\begin{center}
\includegraphics[width=\textwidth]{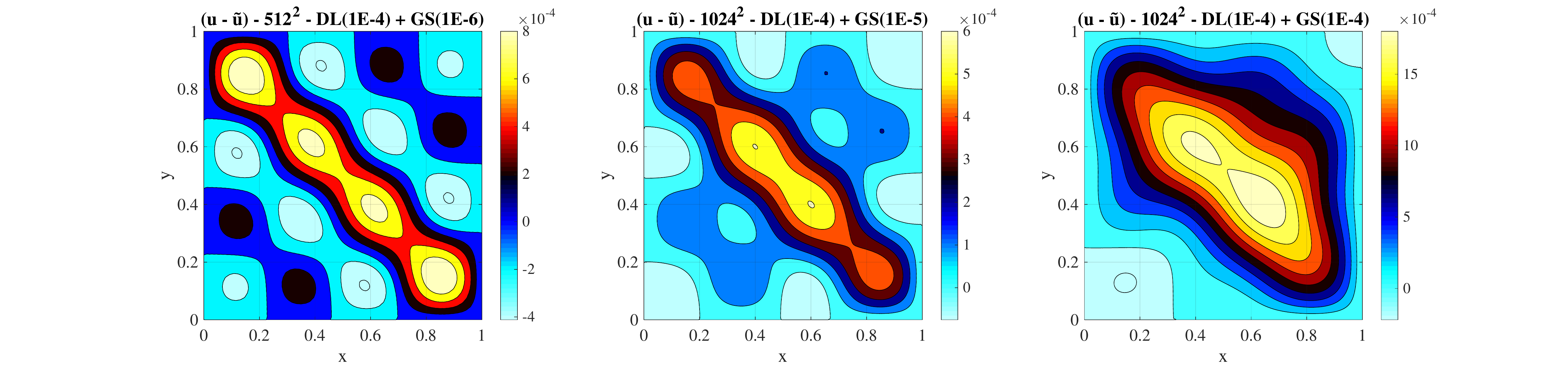}
\end{center}
\caption{Hybrid multigrid solver final error ($u - \tilde{u}$) using three different setups: 1 - 512$\times$512 grid with a 64$\times$64 coarse grid,  $ftol$ equal to 1E-4 and $delta$ equal to 1E-6; 2 and 3 - 1024$\times$1024 grid with a 128$\times$128 coarse grid, $ftol$ equal to 1E-4 and  $\delta$ equal to 1E-5 and 1E-4. }
\label{error1}
\end{figure}

For the first test, we employ a 512$\times$512 grid with a 64$\times$64 coarse grid,  $ftol$ equal to 1E-4 and $\delta$ equal to 1E-6. We then test the hybrid multigrid solver on a 1024$\times$1024 grid with a 128$\times$128 coarse grid, $ftol$ equal to 1E-4 and two values for $\delta$: 1E-5 and 1E-4. Figure \ref{error1} shows a contour plot the error ($u - \tilde{u}$) for these three configurations. The maximum error for the hybrid multigrid solver is of the 1E-4 order and less than the error we obtained after an extensive training of a basic PINN (approximately 1E-3, see the bottom right panel of Figure~\ref{basicPINNresults}). 

\begin{figure}[bt]
\begin{center}
\includegraphics[width=0.9\textwidth]{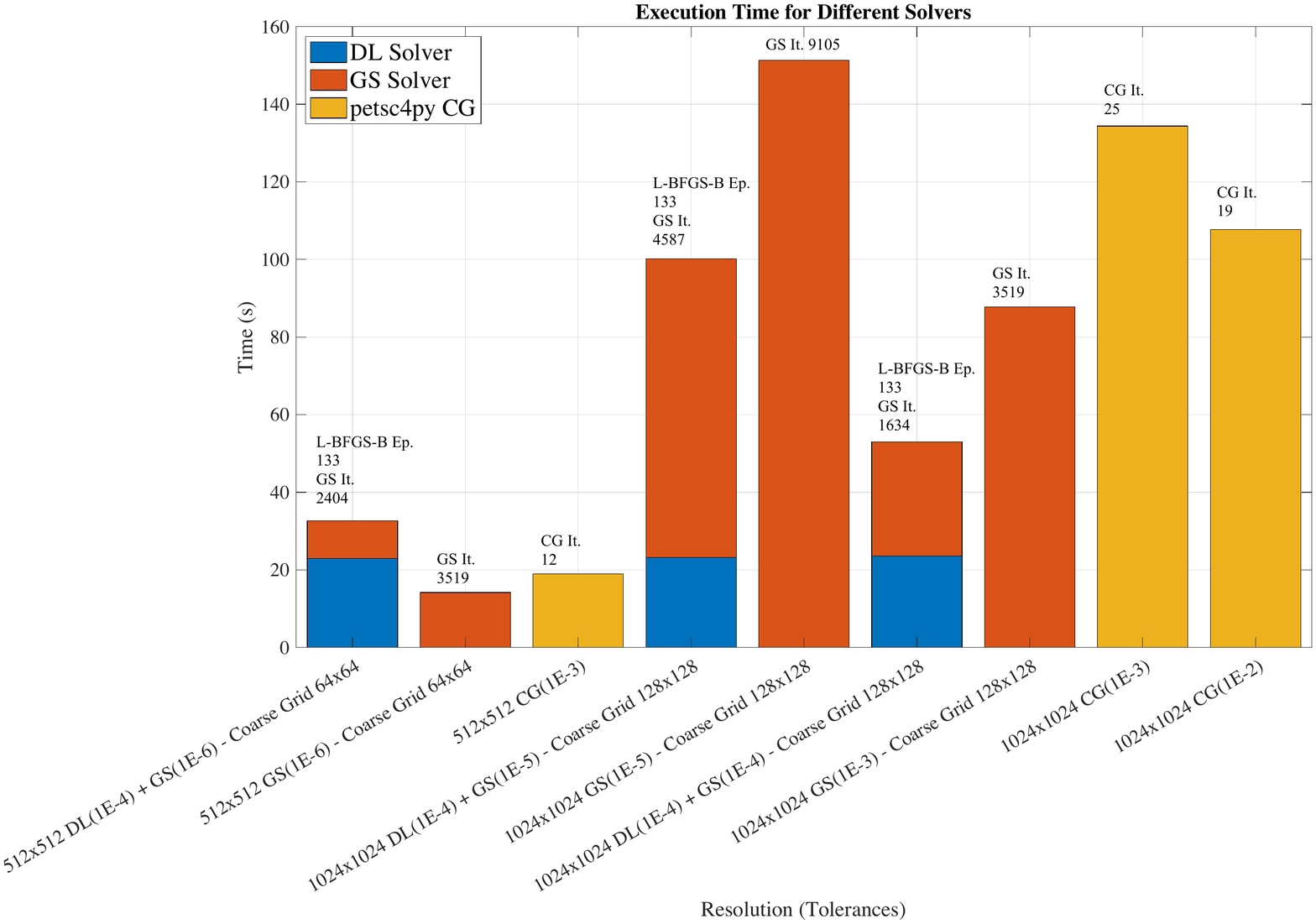}
\end{center}
\caption{Execution time, number of epochs and iterations for the hybrid multigrid DL-GS solver and comparison with the performance of a multigrid using only GS and \texttt{petsc4py} CG varying the resolution, and solver stopping criteria. The hybrid multigrid DL-GS solver is faster for problems using larger coarse grids, e.g. 128$\times$128 coarse grids, than the other approaches.}
\label{performance}
\end{figure}

Once we showed that the hybrid multigrid solver provides more accurate results than the direct PINN usage, we focus on studying the computational performance. The performance tests are carried out on a 2,9 GHz Dual-Core Intel Core i5, 16 GB 2133 MHz LPDDR3 using macOS Catalina 10.15.7. We use Python 3.7.9, \texttt{TensorFlow} 2.4.0, \texttt{SciPy} 1.5.4 and the \texttt{DeepXDE} DSL. The Gauss-Seidel iteration is implemented in \texttt{Cython}~\cite{gorelick2020high} to improve the performance and avoid time-consuming loops in Python.  For comparison, we also solve the problem using only the Gauss-Seidel method to solve the problem on the coarse grid and using the \texttt{petsc4py} CG solver. The \texttt{PETSc }version is 3.14.2 and we use $rtol$ (the relative to the initial residual norm convergence tolerance). We repeat the tests five times and report the arithmetic average of the execution times. We do not report error bars as the standard deviation is less than 5\% of the average value. Figure~\ref{performance} shows the execution time together with number of epochs and iterations for the three different configurations. 

The most important result is that by using an optimized configuration, transfer learning, and integrating DL technologies into traditional approaches, we can now solve the Poisson equation with an acceptable precision with a reduced number of training iterations. This reduction of number of training epochs translates to complete the problem, presented in Section \ref{sec:poisson}, in less than few minutes instead of hours (see Figure~\ref{basicPINNresults}) on the Intel i5 system. While the execution depends on the specific hardware platform and implementation, the number of training epochs and GS iterations on the coarse grid (reported on the top of the histogram bars in Figure~\ref{performance}) are not. Overall, we found that 133 epochs are needed for the L-BFGS-B optimizer to reach an $ftol$ equal to 1E-4.

Figure~\ref{performance} histograms also show the breakdown between the time spent in the DL and Gauss-Seidel solvers used in the multigrid V-cycle. Note that the execution time for the DL solver is approximately the same for calculating the values on the two coarse grids: 64$\times$64 and 128$\times$128. This is because of PINN are \emph{gridless} methods: only the negligible inference computational cost is different. For comparison, we show the performance of the Gauss-Seidel solver for the coarse grid (orange bars) and \texttt{py4petsc} CG solver \texttt{petsc4py} (yellow bars) with different $rtol$ values. When the coarse grid is small, e.g., 64$\times$64, the cost for training the DL solver is higher than using a basic method such Gauss-Seidel: using the Gauss-Seidel method for the coarse grid is faster than using the  DL solver for the coarse grid. However, for larger coarser grids, e.g., 128$\times$128, the hybrid multigrid solver is fastest. For comparison, we present the results obtained running the \texttt{petsc4py} CG with different $rtol$ values. Overall, the performance of the hybrid solver is competitive with state-of-the-art linear solvers. We note that none of the methods and codes have been optimized nor compared at same accuracy (the stopping criteria are defined differently for different solvers), so the performance results provide an indication of potential of the hybrid solver without providing absolute performance values.

\section{Discussion and Conclusion}\label{sec:conclusion}
This paper presented a study to evaluate the potential of emerging new DL technologies to replace or accelerate old traditional approaches when solving the Poisson equation. We show that directly replacing traditional methods with PINNs results in limited accuracy and a long training period. Setting up an appropriate configuration of depth, activation functions, input data set distribution, and leveraging transfer-learning could effectively optimize the PINNs solver. However, PINNs are still far from competing with HPC solvers, such as \texttt{PETSc} CG. In summary, PINNs in the current state cannot yet replace traditional approaches.

However, while the direct usage of PINN in scientific applications is still far from meeting computational performance and accuracy requirements, hybrid strategies integrating PINNs with traditional approaches, such as multigrid and Gauss-Seidel methods, are the most promising option for developing a new class of solvers in scientific applications. We showed the first performance results of such hybrid approaches on the par (and better for large coarse grids) with other state-of-the-art solver implementations, such as \texttt{PETSc}. 

When considering the potential for PINNs of using new emerging heterogeneous hardware, PINNs could benefit from the usage of GPUs that are workforce for DL workloads. It is likely that with the usage of GPUs, the performance of hybrid solvers can outperform state-of-the-art HPC solvers. However, PINN DSL frameworks currently rely on \texttt{SciPy} CPU implementation of the key PINN optimizer, L-BFGS-B, and its GPU implementation is not available in \texttt{SciPy}. The new \texttt{TensorFlow} 2 \texttt{Probability} framework\footnote{\url{https://www.tensorflow.org/probability}} provides a BFGS optimizer that can be used on GPUs. Another interesting research direction is investigating the role and impact of the low and mixed-precision calculations to train the approximator network. The usage of low-precision formats would allow us to use tensorial computational units, such as tensor cores in Nvidia GPUs ~\cite{markidis2018nvidia} and Google TPUs~\cite{jouppi2017datacenter}, boosting the DL training performance.

From the algorithmic point of view, a line of research we would like to pursue is a better and more elegant integration of the DL into traditional solvers. One possibility is to extend the seminal work on discrete PINNs~\cite{raissi2019physics} combining Runge-Kutta solvers and PINN for ODE solutions: a similar approach could be sought to encode information about discretization points into PINN. However, currently, this approach is supervised and requires the availability of simulation data. In addition, the development of specific network architectures for solving specific PDEs is a promising area of research. A limitation of this work is that we considered only fully-connected networks as surrogate network architectures. For solving the Poisson equation and elliptic problems in general, the usage of convolutional networks with large and dilated kernels is likely to provide better performance of fully-connected DL networks to learn non-local relationships a signature of elliptic problems~\cite{lunaaccelerating}.

The major challenge is integrating these new classes of hybrid DL and traditional approaches, developed in Python, into large scientific codes and libraries, often written in Fortran and C/C++. One possibility is to bypass the Python interface of major DL frameworks and use their C++ runtime directly. However, this task is complex. An easier path for the software integration of DL solvers into legacy HPC applications is highly needed. 

Despite all these challenges and difficulties ahead, this paper shows that the integration of new PINNs DL approaches into \emph{old} traditional HPC approaches for scientific applications will play an essential role in the development of next-generation solvers for linear systems arising from differential equations.

\section*{Acknowledgments}
Funding for the work is received from the European Commission H2020 program, Grant Agreement No. 801039 (EPiGRAM-HS).

\bibliographystyle{acm}
\bibliography{PINNscientificComputing}

\end{document}